\documentclass{amsart}
\usepackage[margin=1.5in]{geometry}

\usepackage[all]{xy}
\usepackage{tikz-cd}
\usepackage{enumerate}
\usepackage{amsfonts}
\usepackage{amssymb, amscd, amsmath}
\usepackage{graphicx}
\usepackage{mathtools}
\usepackage{xcolor}
\usepackage{float}
\usepackage[neverdecrease]{paralist}
\usepackage{url}

\newcommand{\R}{\mathbb{R}}

\newcommand{\pf}{\n{\em Proof.}   }

\renewcommand{\t}{\mathfrak{t}}
\newcommand{\bd}{\mbox{bd}}

\newcommand{\inte}{\mbox{int}}

\newcommand{\n}{\noindent}

\newtheorem{teo}{Theorem}[section]
\newtheorem{cor}[teo]{Corollary}
\newtheorem{prop}[teo]{Proposition}
\newtheorem{lema}[teo]{Lemma}

\theoremstyle{definition}

\newtheorem{conj}[teo]{Conjecture}

\title{Flat grazes of convex bodies and local characterization of quadrics}

\author{Jesus Jer\'onimo-Castro }
\author{Luis Montejano   }
\author{Efren Morales-Amaya}

\begin{document}

\begin{abstract}
We prove several local characterizations of quadrics, between them, the local Blaschke's Theorem, and use this result to give some characterizations of the ellipsoid related to the flatness of its grazes.  
\end{abstract}
\maketitle
\section{Introduction}

We think in $\R^n$ as a topological linear vector space over the reals and we extend $\R^n$ to $\mathbb{RP}^n$ by adding the points at infinity. The Grassmann manifold 
of $k$-dimensional subspaces of $\R^n$ is denoted by $\mathcal G(n,k).$  

A line $\ell$ is tangent to a convex body $K \subset \R^{n}$ if the points of $\ell \cap K\not=\emptyset $ lie in the boundary of $K.$  Let us consider a direction $u\in \mathbb{S}^{n-1}$.
The union of all points of the boundary of  $K$,  which lie in some tangent line of $K$ parallel to $u$, will be called  \textit{the shadow boundary } of $K$ in the direction $u$, and will be denoted by $S\partial (K,u).$
That is, 

$$S\partial (K,u):=\{y\in \ell\cap K\mid \ell \mbox{ is a line tangent to } K \mbox{ and parallel to }u\}.$$ 

Abusing notation, if $L$ is a line through the origin in the direction $u$, 
$$S\partial(K,L):=S\partial(K,u).$$ 
\medskip

Let $K\subset\mathbb R^n$ be a convex body, $n\geq 3$. Given a point $x\in\mathbb {RP}^n\setminus K$, define the \emph{graze} of $K$ from $x$:
$$\Sigma(K,x):=\{y\in \ell\cap K\mid \ell \mbox{ is a line tangent  to }K \mbox{ through }x\}.$$
With this notation, the shadow boundary of $K$ in the direction of the line $L$, $S\partial (K,L)$, is the grace $\Sigma (K,x)$, where $x$ is the point at infinity of the line $L$.

\medskip
A set $C\subset \R^n$ is symmetric if and only if $C=-C$.  Moreover, we say that $C$ is centrally symmetric if there is a symmetric translated copy $C+w$ of $C$, that is, if $C$ and $-C$ are translated copies of each other.

 A compact set $C\subset \R^n$  is called a cylinder if it can be represented in the form
$$C= K+L,$$
where $K$ is a $(n-1)$-dimensional convex body in a hyperplane  $\Gamma$,  and $L$ is a line complementary to $\Gamma$.  
The set $K$ is referred to as a base and $L$ as a generatrix of $C$. Note that the generatrix of a cylinder is unique but the base is not. In fact, if $C$ is a cylinder then for every hyperplane   $\Gamma'$ such that the set $K'=C\cap\Gamma'$ is compact, $K'$ is a base of $C$. 

While proving the Banach Isometric Conjecture for $n=3$, Ivanov et. al. (\cite{IMNB} and \cite{IMNK}) proved the following local version of Kakutani`s Theorem:

\begin{teo} [Ivanov]\label{Ivathm}
Let $ K\subset \R^n$ be a convex body containing the origin as an interior point, $n\geq 3$ and let $V$ be an open connected subset of $\mathcal G(n-1,n)$. Suppose that for every $H\in V$, there is $L\in \mathcal G(1,n)$ such that  
$$\bd K\cap H\subset S\partial(K, L).$$
Then there exists a set $\mathcal B\subset \R^n$ such that, for every $H\in V$,

$$K\cap H=\mathcal B\cap H,$$  
and at least one of the following holds:
\begin{enumerate}
\item  there is a line $L_0$ such that $\bd K\cap H\subset S\partial(K, L_0),$  for every $H\in V$ and  
$\mathcal B$ is a cylinder,  or  
\item $\mathcal B=\{v\in \R^n\mid Q(v)\leq 1\}$,  for some non-negative definite quadratic form $Q$ on $\R^n$ .
\end{enumerate}
\end{teo}

In this note, we shall follow the spirit of Ivanov, by proving several local results.  In Section 3, we use these results to study the flat grazes of convex bodies.
In particular, we prove that if $M\subset \inte K$ are $n$-dimensional convex bodies, $n\geq 3$, then $M$ is an
ellipsoid, provided $K$ is a polyhedron and the graces of $M$ from bd$K$ are flat, or that $M$ is an ellipsoid if $M$  
looks supersymmetric from bd$K$.  

\section{Local Results}

We start with the following localized result.

\begin{lema}\label{projection}
Let $K\subset \R^n$ be a smooth convex body, and let $U$ be an open subset of $ \mathcal G(1,n)$, $n\geq 3$.
Suppose that for every $L\in U$, there is $H\in \mathcal G(n-1,n)$ such that
$$\mbox{bd}K \cap H' = S\partial(K, L),$$ 
where $H'$ is parallel to $H$ and also the orthogonal projection of $K$ onto $L^\perp$ is centrally symmetric. 
Then there exists a set $\mathcal B\subset \R^n$ such that, for every $L\in U$,

$$S\partial (K, L) = S\partial(\mathcal B,L)$$  
and  $\mathcal B=\{v\in \R^n\mid Q(v)\leq 1\}$,  for some non-negative definite quadratic form $Q$ on $\R^n$.
\end{lema}

\pf  For the proof of Lemma \ref{projection} we shall use Theorem \ref{Ivathm}. 
Let us first assume that $K$ is smooth.  Define $\Phi: U\to  \mathcal G(n-1,n)$ as follows:
for every $L$ en $U$,  there is $H\in \mathcal G(n-1,n)$ such that  
$$\mbox{bd}K\cap H' = S\partial(K, L),$$ 
where $H'$ is parallel to $H$.  Define  $\Phi(L)= H$.  The function $\Phi$ is well defined,  moreover, $\Phi$ is injective because we assumed that $K$ is smooth.
Consequently $\Phi$ is an open map and therefore $\Phi(U)=V$ is open in $\mathcal G(n-1,n)$.   

For every line $L\in U$,  let $c_L$ be the center of the orthogonal projection $K|L^\perp$,  and let $\ell_L$ be the parallel line to $L$ through $c_L$.  That is, 
$$\ell_L=c_L+L.$$
  We claim that there is $x_0\in K,$ such that:
$$\bigcap_{L\in U} \ell_L \text{ is a single point } \{x_0\}.$$

Let $L_1$ and $L_2$ be two different lines through the origin, such that $L_i\in U, i=1,2.$   The orthogonal projection of $K$ onto $(L_1^\perp\cap L_2^\perp)$ is centrally symmetric, because by hypothesis $K|L_i^\perp$ is centrally symmetric, $i=1,2$; let $a$ be its center of symmetry and
let $\Gamma$ be the orthogonal plane to $L_1^\perp\cap L_2^\perp$ passing through $a$. That is 
$$\Gamma= a+(L_1+L_2).$$ 
The lines $\ell_{L_1}$ and $\ell_{L_2}$ are both in $\Gamma$ because the symmetry centers $c_{L_1}$ and $c_{L_2}$ are orthogonally projected to $a$. Therefore, since 
$\ell_{L_1}$ and $\ell_{L_2}$ are not parallel, they intersect in a point $\{x_0\}$. 
Consider a third line $L_3\in U$ through the origin, linearly independent to $L_1$ and $L_2$; it exists because $n\geq 3$. By the above argument, $\ell_{L_3}$ intersects both $\ell_{L_1}$ and $\ell_{L_2}$. But by linear independence, $\ell_{L_3}$ intersects $\Gamma$ in at most one point, so we must have that $x_0\in \ell_{L_3}$.
Finally, for the general line $L\in U$, $L$ is linearly independent to at least one pair of the lines $L_1, L_2, L_3$. So that the preceding argument yields that $x_0\in \ell_{L}$. This proves that  $\bigcap_{L\in U} \ell_L \text{ is a single point } \{x_0\},$ as we claimed.

Without loss of generality we may assume that $x_0$ is the origin. Hence, by the above, for every $L\in U$, $K|L^\perp$ is centrally symmetry with center at the origin. 
Furthermore, since by hypothesis, $\mbox{bd}K\cap H' = S\partial(K, L),$  where $H'$ is parallel to $H$, we have that $K\cap H'$ is centrally symmetric with with center at $b_L$ and the orthogonal projection onto $L^\perp$ sends $b_L$ to the origin.  

Our next purpose is to prove that $K\cap H' $ contains the origin. By Theorem \ref{Ivathm}, this will be enough to prove our Lemma, because for every $H\in V$,   there is $L\in \mathcal G(1,n)$ such that  
$\mbox{bd}K\cap H=S\partial(K, L),$ where we know that $V$ is open in $\mathcal G(n-1, n)$. Note that by hypothesis $\mathcal B$ is not a cylinder.

Let $[p,q]$ be a diametral chord of $K\cap H'$.  We shall prove that 
$\frac{p+q}{2}$, the midpoint of  $[p,q]$ is the origin, thus proving that  the origin lies in $K\cap H'$. 
 By definition there are support hyperplanes $p+\Delta$ and $q+\Delta$ of $K$ at $p$ and $q$, where $\Delta$ is a hyperplane through the origin containing $L$. 
 Let $U'=\{L\in \mathcal G(1,n)\mid L\subset \Delta\}$. So we may assume that $U'$ is an open subset of the set of lines through the origin in $\Delta$.
   If $\pi:\R^n\to L^\perp$ is the orthogonal projection, where $L\in U'$, then $\pi([p,q])$ is a diametral chord of $\pi(K)$, whose midpoint is the origin because $\pi(K)$ is centrally symmetric with center at the origin. This implies that $\pi(p)=-\pi(q)$. Therefore, 
$\pi(\frac{p+q}{2})=0$, for every orthogonal projection  $\pi:\R^n\to L^\perp$, with $L\in U'$, which may be considered as an open subset of the set of lines through the origin in $\Delta$.  
 Consequently, $\frac{p+q}{2}$ must be the origin. \qed

\medskip

Next, we state the  local version of Blaschke's  Theorem. See Section 2.12 of \cite{MMO}.

\begin{teo} \label{Blthm}
Let $K\subset \R^n$ be a  smooth convex body and let $U$ be an open subset of $ \mathcal G(1,n)$, $n\geq 3$.
Suppose that for every $L\in U$, there is hyperplane $H$  such that  
$$\mbox{bd}K\cap H = S\partial(K, L),$$  
Then there exists a set $\mathcal B\subset \R^n$ such that, for every $L\in U$,

$$S\partial (K, L) = S\partial(\mathcal B,L)$$  
and  $\mathcal B=\{v\in \R^n\mid Q(v)\leq 1\}$,  for some non-negative definite quadratic form $Q$ on $\R^n$.
\end{teo}
\pf As in the first part of the proof of Lemma \ref{projection}, define $\Phi: U\to  \mathcal G(n-1,n)$ in such a way that $\Phi$ is an open map.

We shall prove our theorem by induction. First suppose $n=3$. Let $L_0\in U\subset \mathcal G(1, 3)$, our main step is to prove that there is a neighborhood 
$W$ of $L_0$ in $U$ and a quadric $\mathcal B\subset \R^3$ such that, for every $L\in W$,
$S\partial (K, L) =S\partial(\mathcal B,L).$
Since two quadrics that coincide in an open set are the same, this concludes the proof of the theorem. For the proof of the main step, it will be enough to prove that there is a neighborhood $W$ of $L_0$ in $U$ and a point $x_0\in $ int$K$, such that for every $L$ in $W$, the shadow boundary  $S\partial(K,L)$ is flat and furthermore, if $\mbox{bd}K\cap H = S\partial(K, L)$, then $x_0 \in K\cap H $.  Without loss of generality, we may assume that $x_0$ is the origin and consequently, by Theorem \ref{Ivathm}, for $V=\Phi(W)$, since $K$ is strictly convex, 
 there exists a set 
$\mathcal B\subset \R^{3}$ such that, for every $L\in W$,
$$S\partial ( K, L) = S\partial(\mathcal B,L)$$  
and  $\mathcal B=\{v\in \R^n\mid Q(v)\leq 1\}$,  for some non-negative definite quadratic form $Q$ on $\R^{3}$.

For the proof of the existence of $W$, consider the closed set $\mathcal A=\{\ell \in \mathcal G(1,3)\mid \ell\subset \Phi(L_0)\}$. 
 If $\Omega$ is a small open neighborhood of $\mathcal A$ in $\mathcal G(1,3)$, then $\tilde \Omega\subset \mathcal G(2,3)$ denotes the open subset of $G(2,3)$ consisting of all planes generated by a pair of different lines of $\Omega$. 
 Let $\Omega$ be so small that $\tilde\Omega\subset \Phi(U).$
Let's  prove that there is a point $x_0\in \R^{3}$ with the property that if $[p,q]$ is a diametral chord of $K$ and 
the line through the origin parallel to Aff$[p,q]$ is in $\Omega$ then $x_0\in [p,q]$.  Consider two diametral chords $[a,b]\not= [a',b']$ of $K$, where the line $\ell$ through the origin, parallel to Aff$[a,b]$ is in $\Omega$
and the line $\ell'$ through the origin parallel to Aff$[a',b']$ is in $\Omega$.  We know that  the plane $H\in \mathcal G(2,3)$ generated by $\ell$ and $\ell'$ lies in $\Phi(U)$
 Consequently, there is a direction $\mathcal L$ for which 
$$S\partial(K,\mathcal L)= H'\cap \mbox{bd} K,$$
with $H'$ parallel to $H$.  Since $[a,b]$ and $[a',b']$ are parallel to $H'$ and both diametral chords, then both are contained in the plane $H'$ and therefore they have non empty intersection.  Denote by $x_0=[a, b]\cap[a', b']$.  Let now $[a'',b'']$ be a diametral chord of $K$ such that the line 
$\ell''$ through the origin, parallel to Aff$[a'',b'']$ is in $\Omega$ and the lines $\{\ell,\ell',\ell''\}$ are linearly independent. By the above, $x_0\in[a'',b'']$.  Finally, let $[p,q]$ be a diametral chord of $K$ such that the line 
$L$ through the origin, parallel to Aff$[p,q]$ is in $\Omega$. Since $L$ is linearly independent to at least one pair of $\{\ell,\ell',\ell''\}$, then $x_0\in[p,q]$. 
By continuity of $\Phi$, let $W$ be a neighborhood of $L_0$ in $U$, such that $\Phi(W)\subset \tilde\Omega\subset\phi(U)$. Hence, if $L\in W$,  there is $H'$ such that $S\partial (K,L)=\mbox{bd}K\cap H'$ with $H'$ parallel to $H\in \tilde\Omega$. The fact that every diametral chord of $K$ parallel to a line in $\Omega$ pass through $x_0$, implies 
$H'$ also pass through $x_0$. This conclude the proof of the main step and therefore the proof of the theorem when  $n=3$.

Suppose now the Theorem is true for $n-1$, we shall prove it for $n$.  
By hypothesis there is a open set $U\subset \mathcal G(1,n)$ such that for every $L\in U$, there is $H\in \mathcal G(n-1,n)$ such that 
$$S\partial( K, L)=H'\cap \mbox{bd} K.$$
where $H'$ is parallel to $H$.  We shall show that $H'\cap  K$ is an ellipsoid.   For that purpose, let  $\Gamma$ by a hyperplane, parallel to $L$ that intersects the interior of  $H'\cap  K$.
We shall use our induction hypothesis for $\Gamma \cap K$ in $\Gamma$.  Clearly  the shadow boundary of $\Gamma \cap  K$ in $\Gamma$, 
$$S\partial (\Gamma \cap  K,L)=S\partial( K, L)\cap \Gamma=(H'\cap \Gamma)\cap \mbox{bd} K= 
(H'\cap \Gamma) \cap \mbox{bd} (\Gamma \cap  K).$$
In other words, the shadow boundary of $\Gamma \cap K$ in the direction of $L$ is flat, but the same is true for every $L'\in U$, parallel to $\Gamma $.  By induction
$(H'\cap \Gamma)\cap K$ is an ellipsoid. Since this holds for every hyperplane $\Gamma$ parallel to $L$ that intersects the interior of  
$H'\cap K$, we have that for every $L\in U$,
  $$S\partial( K, L)=H'\cap \mbox{bd}K \mbox{  is an ellipsoid  }.$$

Remember that $\Phi (U)$ is open in $\mathcal G(n-1,n)$ and consider, for every $H\in \Phi (U)$, the orthogonal projection of $K$ onto $L^\perp$, where
 $S\partial( K, L)=H'\cap \mbox{bd} K$ and  $H'$ is parallel to $H$.  Then $ K|L^\perp=(H'\cap K)|L^\perp$ and hence, for every $H\in \Phi(U)$, 
$ K|L^\perp$ is an ellipsoid.  By Lemma \ref{projection}, there exists a set $\mathcal B\subset \R^n$ such that, for every $H\in \Phi (U)$,

$$ H'\cap  \mbox{bd}K= S\partial( K, L)=S\partial(\mathcal B, L)=H'\cap \mbox{bd}\mathcal B$$  
and  $\mathcal B=\{v\in \R^n\mid Q(v)\leq 1\}$,  for some non-negative definite quadratic form $Q$ on $\R^n$. 
\qed

\bigskip

In 1977 Burton \cite{Bur2} proved that a convex body all whose sections are centrally symmetric must be an ellipsoid. Next, we have a local version of this result. Denote by $\mathcal G^0(n-1,n)$ the Grassmannian space of all hyperplanes
(not necessarily through the origin)  of $\R^n$.

\begin{teo} \label{thmBur}
Let $K\subset \R^{n}$ be a  convex body containing the origin as an interior point, $n\geq 3$. Let $V$ be an open subset of 
$\mathcal G^0(n-1,n)$ such that for every 
$H\in V$, the section $H\cap K$ is a non-empty centrally symmetric convex body. Moreover, suppose $\bigcup_{H\in V} (H\cap K)$ is symmetric.  Then there exists a set $\mathcal B \subset \R^n$ such that, for every $H\in V$,
through the origin,
$$H\cap K=H\cap \mathcal B,$$ 
where $\mathcal B=\{v\in \R^n\mid Q(v)\leq 1\}$,  for some non-negative definite quadratic form $Q$ on $\R^n$.
\end{teo}

\pf  Let us consider first the  
following situation:  let $M\subset \R^n$ be a strictly convex body and let $\Gamma$ by a hyperplane intersecting the interior of $M$. Suppose that 
$$(\Gamma \cap M)+ v= (\Gamma +v)\cap M.$$
for some vector $v\in \R^n$.  Let $B$ be the band between the hyperplanes  $\Gamma$ and $\Gamma + v$ and consider the line  $L=\{tv\mid t\in\R\}$. Then is is clear that
$$M\subset \big((\Gamma\cap M) +L\big)\cup B.$$
Furthermore, $S\partial (M,L)\subset B.$

The next step is to use Theorem \ref{Ivathm}. 
We shall prove that for any hyperplane  $H\in V$
through the origin, the boundary of $H\cap K$ is a shadow boundary of $K$.
For that purpose let  $K_{i}=(H+\frac{\epsilon}{i}u)\cap K$ and $%
K_{-i}=(H-\frac{\epsilon}{i}u)\cap K$ in such a way that $\epsilon$ is so small that $H+\frac{\epsilon}{i}u$ and $
H-\frac{\epsilon}{i}u$ are hyperplanes of $V$. Consequently, by hypothesis $K_i$ and $K_{-i}$ are centrally symmetric. Furthermore, since  $\bigcup_{H\in V} (H\cap K)$ is symmetric, then 
$K_{-i}$ is a translated copy of $K_{i}.$  Assume that  
$$K_{i}=K_{-i}+v_{i}.$$ 
Without loss of generality, we may assume that the vectors $\frac{v_{i}}{%
\parallel v_{i}\parallel }$ converge to a vector $v$. 

By the claim in the beginning of the proof,
$$K\subset \big(K+L_i\big)\cup (B_i\cap M),$$
where $L_i=\{tv_i\mid t\in\R\}$ and $B_i$ is the band between the hyperplanes $(H+\frac{1}{i}u)$ and $(H-\frac{1}{i}u).$
Furthermore, $S\partial (K,L_i)\subset B_i.$

Consequently, in the limit when $v_i\to v,$  $K\subset (H\cap \mbox{bd }K)+L$ and $S\partial(K,L)\subset H,$
where $L=\{tv\mid t\in\R\}$. Therefore $H\cap \mbox{bd} K=S\partial(K,L)$.  Since the same holds for every hyperplane $H\in V$ through the origin, by Theorem \ref{Ivathm}, 
there exists a set $\mathcal B \subset \R^n$ such that, for every $H\in V$,
through the origin
$$H\cap K=H\cap \mathcal B$$ 
Note that   $\mathcal B$ is not a cylinder, hence $\mathcal B=\{v\in \R^n\mid Q(v)\leq 1\}$,  for some non-negative definite quadratic form $Q$ on $\R^n$.
\qed
\bigskip

T. Kubota \cite{Ku2} proved that a convex body 
is an ellipsoid provided any two of its sections by parallel hyperplanes are direct homothetic. We present here a local version of this result

\bigskip

\begin{teo}
Let $K\subset \R^n$ be a smooth convex body and let $U\subset \mathcal G(1,n)$ be an open set.  
Suppose that any two parallel sections $H_1\cap K$ and $H_2\cap K$ of $K$ 
are directly homothetic, whenever there is $L\in U$ parallel to $H_1$ and $H_2$.
Then, there exists a set $\mathcal B \subset \R^n$ such that for every $L\in U$, 
$$S\partial ( K, L)= S\partial (\mathcal B, L),$$
where $\mathcal B=\{v\in \R^n\mid Q(v)\leq 1\}$,  for some non-negative definite quadratic form $Q$ on $\R^n$.
\end{teo}
\pf Suppose first $n=3$.  By Theorem \ref{Blthm}, it will be enough to prove that for every $L\in U$,  $S\partial ( K, L)$ lies in a plane. 
  Let $S$ be the set of all unit vectors orthogonal to $L,$ and let $v\in S$. Suppose
$H$ is a plane orthogonal to $v$, which cuts the interior of $K.$ As we assume
that $K$ is strictly convex, then $H\cap S\partial
(K,L)$ has exactly two points, say $a$ and $b.$ If this is the case, we say that the chord
$ab$ of $S\partial (K,L)$  is a  $v$-chord of
$S\partial (K,L).$ We know that the support lines of $H\cap K$ parallel to $L$ pass through $a$ and $b$.  Therefore the chord $ab$
is a diametral chord of $H\cap K$. 
If $H_{1}$ and $H_{2}$  are two planes orthogonal to $v$,  the $v$-chords $H_{1}\cap S\partial (K,u)$ and $H_{2}\cap S\partial (K,u)$ of  
$S\partial (K,u)$ are parallel, since $H_{i}\cap S\partial (K,u)$ is  a chord of $H_{i}\cap K$ in the direction $u,$ $%
i=1,2,$  and  $H_{1}\cap K$  is directly homothetic to $H_{2}\cap K$.

This proves that for any $v\in S$ there is a line $L_v$ through the origin, such that all the   $v$-chords of $%
S\partial (K,u)$ are parallel to $L_v$. This wil be enough
to claim that  $S\partial (K,L)$ lies on a plane.  To prove this, our next step is to show that there is a plane $T$ through the origin, such that for every $v\in S,$

$$L_v= T\cap v^\perp.$$

 Let $v_{1}, v_{2}$ in $S$ and let $T$ be the plane generated by  $L_{v_{1}}$ and  $L_{v_{2}}.$
Let us consider now the inverse Gauss map $\gamma: \mathbb{S}^2\to \mbox{bd}K$, which assigns to every unit vector $w$ the unique point  $\gamma(w)$ of  $\mbox{bd}K$
with the property that $w$ is a normal vector to $\mbox{bd}K$ at $\gamma(w)$. Note that 
$S\partial (K,L)=\{\gamma(v)\mid v\in S\}.$ Take  $c\in S\partial (K,L)$ between 
$\gamma(v_{1})$ and  $\gamma(v_{2})$ and let $ca$ be the  $v_{1}$-chord of  $S\partial (K,L)$ and
 $cb$  be the  $v_{2}$-chord of $S\partial (K,L).$
Let $H'$ be the plane containing the chord  $ab$ of $S\partial (K,L)$ and parallel to $L$, and suppose that $w$ is orthogonal to $H'$. By definition $w\in S$,
$ab$ is a $w$-chord of $S\partial (K,u)$, and at least for this unit vector $w\in S$, the line $L_{w}= T\cap {w}^\perp$.

Clearly, $w$ depends continuously on the point $c$ while moving along $S\partial (K,L).$
 Note that when $c$ approximates $\gamma(v_{2})$,   $a$ also approximates $\gamma(v_{2})$. Hence
$w$ approximates $v_{1 }$. Similarly, when $c$
approximates $\gamma(v_{1})$ also $b$ does, therefore $w$ approximates $v_{2}$.
Consequently, while $c$ travels between $\gamma(v_{2})$ and $\gamma(v_{1})$, $w$
travels continuously between $v_{1}$ and $v_{2}$. It is easy to verify  that while $c$ travels through  $S\partial (K,L)$, the corresponding $w$ takes all the directions of $S$.
This shows that for every unit vector
$w\in S$, $L_w= T\cap w^\perp$.

Finally, this gives us enough information to prove that $S\partial (K,L)$ lies on a plane because if $p\in S\partial (K,L)$ is a fixed point, then

$$S\partial (K,L)\setminus\{p\}=\{ q\in \mbox{bd}K \mid [p,q] \mbox{ is a $v$-chord}, v\in S\setminus\{\gamma^{-1}(p)\}\}.$$
Given that $[p,q]\subset p+L_v \subset p + T$, we have that $S\partial (K,L)$  $=(p+T)\cap \mbox{bd} K$, as we wanted.

The proof of case $n>3$ is straightforward but using the notion of $k$-shadow boundary instead of the notion of shadow boundary.
 \qed

\bigskip

Peter R. Goodey \cite{Goo} proved that if for every translated copy $K'\not=K$ a a convex body, the intersection 
$\mbox{bd} K\cap \mbox{bd} K'$ is flat, then $K$ is an ellipsoid. Here is the corresponding local version.

\bigskip
\begin{teo}

Let $K\subset \R^n$ be a smooth convex body, $n\geq 3$, and let $U\subset \mathcal G(1,n)$ be an open set and let $\Omega_U=\{v\in \R^n \mid v \mbox{ is a nonzero vector in some of the lines of } U$.   Suppose that 
$$\bd K\cap \bd K'$$
is contained in a hyperplane whenever $K'=K+v$ and $v\in \Omega_U$
Then, there exists a set $\mathcal B \subset \R^n$ such that for every $L\in U$, 
$$S\partial ( K, L)= S\partial (\mathcal B, L),$$
where $\mathcal B=\{v\in \R^n\mid Q(v)\leq 1\}$,  for some non-negative definite quadratic form $Q$ on $\R^n$.
\end{teo}
\pf  Suppose $\mbox{bd} K\cap \mbox{bd}K'$ is flat, where $K$ is a convex body and $K'$ is a translated copy of  $K$ such that  $\mbox{bd}K'\cap$ int$K$ is not empty. Then, there is a hyperplane $H$ such that $H\cap\mbox{bd} K=H\cap\mbox{bd} K'=\mbox{bd}K\cap \mbox{bd} K'$. 
If there is a hyperplane $H$ such that $\mbox{bd} K\cap \mbox{bd}K'\subset H$, then clearly, $\mbox{bd }K\cap \mbox{bd}K' \subset \mbox{bd }K\cap H$.  Suppose that  $\mbox{bd }K\cap \mbox{bd}K'$  is properly contained in  
 $\mbox{bd} K\cap H$. This implies that  $\mbox{bd} K \setminus \mbox{bd} K'$ is not connected, but this is a contradiction to the fact that  
 $\mbox{bd}K'\cap$ int$K$ is not empty. Similarly, $\mbox{bd} K\cap \mbox{bd}K'= H\cap \mbox{bd}K'$. 
 
 Our next purpose is to prove that if $\mbox{bd} K\cap \mbox{bd}(K-tv)$ is flat for every $|t|<\epsilon$, then $S\partial(K, v)$ is flat.  Let $\epsilon>0$ be so small that $\mbox{bd}(K-tv)\cap$ int$K$ is not empty. By hypothesis,  for every $|t|<\epsilon$, there is a hyperplane $H_t$ such that 
 $H_t\cap\mbox{bd} K=H_t\cap\mbox{bd}( K-tv)=\mbox{bd}K\cap \mbox{bd}(K-tv)$. Therefore, 
 $$(H_t\cap K)+tv=(H_t+tv)\cap K.$$
By the claim at the beginning of the proof of Theorem \ref{thmBur},  $K\subset \big((H_t\cap K)+L)\cup B_t$, where $L$ is the line $\{tv\mid t\in\R\}$, and $B_t$ is the band between the parallel hyperplanes $H_t$ and $H_t+tv$. Furthermore $S\partial (K,v)\subset B_t$.  Choose a sequence $\t_i\}$ converging to $0$, such that 
$\lim (H_{t_i})=\lim (B_{t_i})=H$.  Consequently, $S\partial (K,L)= S\partial (K,v)\subset H,$ as we wished. 

Finally, since this is true for every $L\in U$, By Theorem \ref{Blthm}, there exists a set $\mathcal B \subset \R^n$ such that for every $L\in U$, 
$$S\partial ( K, L)= S\partial (\mathcal B, L),$$
where $\mathcal B=\{v\in \R^n\mid Q(v)\leq 1\}$,  for some non-negative definite quadratic form $Q$ on $\R^n$. \qed 

\bigskip
\bigskip

\section{Flat graces of convex sets}

Let $M\subset\mathbb R^n$ be a convex body, $n\geq 3$. Given a point $x\in\mathbb R^n\setminus M$, 
Remember that the \emph{graze} of $K$ from $x$:
$$\Sigma(M,x):=\{y\in L\cap K\mid L \mbox{ is a support line to }M \mbox{ through }x\}.$$

Suppose $H$ is a hyperplane not intersecting the convex body $M$  and suppose $\Sigma(M,x)$ is flat for every $x\in H$. In 1959, Marchaud \cite{March}, proved that if this is the case $M$ must be an ellipsoid.  In 1977, Burton \cite{Bur2} showed that if for some $\delta >0$, the set of points in $\R^n\setminus M$ whose distance from $M$ is less than $\delta$, have flat grazes then $M$ is an ellipsoid.

The following conjecture, which is still open, seems to be very natural.  

\begin{conj}\label{conjBG}

Let $M, K$ be convex bodies in $\R^n$ and suppose $M\subset$ int$K, n\geq3 $. If the graze $\Sigma(M,x)$ is flat  for every $x\in \mbox{bd} K$, then $M$ is an ellipsoid.
\end{conj}

\bigskip

\subsection{Flat grazes via the local  Blaschke's Theorem}
The purpose of this section is to prove that Conjecture \ref{conjBG} is true under certain additional hypotheses, in particular if $K$ is a polyhedron.  As an immediate corollary we obtain the smooth version of Burton's  Theorem \cite{Bur2} and also the following variation of Conjecture \ref{conjBG},  \emph{Suppose $M$ is smooth and $M\subset \inte K$, 
if for some $\delta >0$, the set of points whose distance from $\bd K$ is less than $\delta$, have flat grazes then $M$ is an ellipsoid.}

\bigskip

Our first result is

\begin{teo}\label{poly}
Let $M, P$ be convex bodies in $\R^n$ and suppose $M$ is smooth and $M\subset \inte P$, where $P$ is a polyhedron, $n\geq 3$.  If the graze $\Sigma(M,x)$ is flat  for every $x\in \bd P$, then $M$ is an ellipsoid.
\end{teo}
\pf  Let $x\in \bd P$  and suppose $x$ is in the interior of a maximal face $\sigma$ of the boundary of $P$.  Suppose $\Gamma$ is a hyperplane with the property that $\Gamma\cap P=\sigma$.  Since $\Gamma$ is a support hyperplane of $P$ and 
$M\subset$ int$P$, then $\Gamma\cap M=\emptyset$. Let $\Phi$ be a projective isomorphism that sends $\Gamma$ to the hyperplane at infinity and denote by $\mathcal M$ the image of $M$ under $\Phi$. Let $L_0$ a line through the origin such that 
the point at infinity of $L_0$ is $\Phi(x)$.  Hence, by hypothesis, there is an open neighborhood $U$ of $L_0$ in $\mathcal G(1,n)$
such that  that for every $L\in U$, there is $H\in \mathcal G(n-1,n)$ such that  
$$\mbox{bd}\mathcal M\cap H' = S\partial(\mathcal M, L),$$ 
where $H'$ is parallel to $H$.    By Theorem \ref{Blthm},
there exists a set $\mathcal B'\subset \R^n$ such that, for every $L\in U$, $S\partial (\mathcal M, L) = S\partial(\mathcal B',L)$  
and  $\mathcal B'=\{v\in \R^n\mid Q(v)\leq 1\}$,  for some non-negative definite quadratic form $Q$ on $\R^n$. Consequently, $S\partial (\mathcal M, L)$ is an ellipsoid and therefore the graze $\Sigma (M,x)$ is also an ellipsoid.  
Suppose now $x_i\in $ int$\sigma$ is a sequence of points in the interior of a maximal face of $P$ such that  $\lim x_i=y$.  By the above $\Sigma(M,x_i)$ is an $(n-1)$-ellipsoid for every $i\geq 1$ and $\lim \Sigma(M,x_i)=\Sigma(M,y)$.
Therefore $\Sigma(M,y)$ is also  an $(n-1)$-ellipsoid. This implies that for every $x\in P$, the cone of $M$ from {x} is an ellipsoidal cone and hence,
by the Bianchi-Gruber Theorem 2 in \cite{BG}, $M$ is an ellipsoid as we wished  \qed

\bigskip

Next we have several immediate corollaries to Theorem \ref{poly}

\begin{cor}
Let $M, K$ be convex bodies in $\R^n$ and suppose $M$ is smooth and  $M\subset$ int$K, n\geq3 $.  If for some $\delta >0$, the set of points whose distance from bd$K$ is less than $\delta$, have flat grazes, then $M$ is an ellipsoid.
\end{cor}
\pf  Approximate $K$ with a polyhedron. \qed

\medskip

\begin{cor}[Burton]
Let $M$ be a smooth convex body in $\R^n, n\geq3 $.  If for some $\delta >0$, the set of points in $\R^n\setminus M$ whose distance from $M$ is less than $\delta$, have flat grazes then $M$ is an ellipsoid.
\end{cor}
\pf  Approximate $M$  from outside with a polyhedron. \qed

\medskip

Below we state several characterizations of the ellipsoid related to the fact that from outside a convex body there are many points that are its poles.

We complete  $\mathbb{R}^{n}$, by adding the points at infinity,  to obtain  projective
$n$-dimensional space $\mathbb{RP}^{n}.$
If $x\in \mathbb{RP}^{n}\setminus \mbox{bd} M,$ 
$x$ is a \textit{pole } of \textit{\ }$M$ if there exists a hyperplane $H$ of $%
\mathbb{RP}^{n}$ with the property that for every line  $L$ trough $x$ intersecting the interior of $M$ and such that $%
\mbox{bd}M\cap L=\left\{ a,b\right\} $, the cross ratio of 
 $A,B,x$ and $L\cap H$ is minus one. That is:

\begin{center}
$\left[ a,b;x,p\right] =\frac{{\overline{ax}}^.{\overline{pb}}}{{\overline{xb}}^.{\overline{ap}}}=-1.$
\end{center}

\noindent where $p=L\cap H$ and  $\overline{ax}$, $\overline{pb}$, $\overline{xb}$ and $\overline{ap}$ denotes the signed length of the directed chords $[a,x]$, $[p,b]$, $[x,b]$ and $[a,p]$, respectively, in the directed line $L$.  If this is so, we say that $H$ is a \textit{polar hyperplane}
of $M$ and also that $H$ is \textit{the polar of the pole} $x$.  See \cite{MM1}.

Remember that projectivities preserve the cross ratio and that 
$\left[ a,b;x,p\right] =-1,$ with  $p=\infty $ if and
only if $x$ is the midpoint of the interval $[a,b]$.  Consequently,
if $x\in \mathbb{RP}^{n}\setminus M$ is a pole of $M$, then its polar $H$ is 
a projective hyperplane of affine symmetry of $M$. In this case, 
$$H\cap \mbox{bd} M= \Sigma (M,x).$$  
As a consequence, we have the following two corollaries

\medskip
\begin{cor}
Let $M, P$ be convex bodies in $\R^n$ and suppose $M$ is smooth and $M\subset$ int$P$, where $P$ is a polyhedron, $n\geq 3$.  If every point $x\in \mbox{bd} P$ is a pole of $M$, then $M$ is an ellipsoid.
\end{cor}

\medskip

\begin{cor}
Let $M, K$ be convex bodies in $\R^n$ and suppose $M$ is smooth and $M\subset$ int$K, n\geq3 $.  If for some $\delta >0$, the set of points whose distance from bd$K$ is less than $\delta$, are poles of $M$, then $M$ is an ellipsoid.
\end{cor}

\bigskip

\subsection{ Convex bodies that look supersymmetric}
Let $M\subset \R^n$ by a symmetric convex body and let $x\notin M$. We say that $M$ looks supersymmetric from $x$ if the grace $\Sigma(M,x)$ is planar and centrally symmetric whose center lies between the origin and $x$.
An ellipsoid is supersymmetric from any exterior point from which it is viewed.

Note that if $M$ looks supersymmetric from $x$, then according with Bianchi and Gruber \cite{BG},  $M$ looks symmetric from $x$, because the line through $x$ and the origin is a pole of the cone of $M$,  $C_x(M)$.  Bianchi and Gruber \cite{BG} conjectured that a convex body $M$ that looks symmetric from 
a surface surrounding $M$ should be an ellipsoid.  Indeed, Gruber and Odor \cite{GO} proved that this is the case if $M$ is of class $\mathcal C^4$ and looks symmetric from every point sufficiently close to $M$:

Note that our notion of supersymmetry is related, by duality, with the False Pole Theorem (we recommend the interested reader to go to \cite{larmanmorales} )

The purpose of this section is to prove the following result:

\begin{teo}\label{super}
  Let $M\subset \R^n$ by a symmetric convex body and let $\mathcal E$ be an ellipsoid centered at the origin such that $\sqrt{2} M\subset \mbox{int}\mathcal E$, $n\geq 3$. If $M$ looks supersymmetric from every point of the boundary of $\mathcal E$, Then $M$ is an ellipsoid.
\end{teo}

Let $M\subset \R^n$ be a convex body and suppose $x\in \R^n\setminus M$. Define the \emph{double cone}:
$$S(M,x):=\{y\in L\mid L \mbox{ is a support line to }M \mbox{ through }x\},$$
in such a way that $S(M, x)\cap M=\Sigma(M,x)$.

The \emph{saturn-ring}  $\Omega(M,x)$ of $M$ with respect to $x$ is:
$$\Omega(M,x)=: S(M,x)\cap S(M-x).$$

Let $M$ and $K$ be two symmetric convex bodies in $\R^n$, 
with $M \subset  K$. We say that the body $K$ is \emph{ almost free with respect to} $M$  if for each $x \in \mbox{bd} K$, the saturn ring of $M$ with respect $x$ is contained in $K$, that is:
$$\Omega(M,x)\subset K.$$

Suppose $\Sigma(M,x)$ is a planar centrally symmetric  whose center lies between the origin and $x$.
Let us prove that $\Omega(M,x)$ is flat. First of all note that $\Omega(M,x)$ is symmetric. Let $y \in \Omega(M,x)$ and let $P$ be the plane $\mbox{aff }\{x,-x, y\}$.   Since $L(y,x)$ is a support line of $P\cap M$, let $a=L(y,x)\cap P\cap M$.
Similarly, $L(y,-x)$ is a support line of $P\cap M$, let $b=L(y,-x)\cap P\cap M$. By symmetry,  $L(-y,x)$ is a support line of $P\cap M$, let $-b=L(y,-x)\cap P\cap M$ and moreover, $L(-y,-x)$ is a support line of $P\cap M$, let $-a=L(-y,-x)\cap P\cap M$.
Therefore, $-y\in \Omega(M,x)$. Of course $\{a, b\}=P\cap \Sigma(M,x).$  
Since $[0,x]\cap [a,b]$ is the midpoint of $[a,b]$, then $[-y,y]$ is parallel to $[a,b]$.
If the graze $\Sigma(M,x)$ is contained in the hyperplane $\Delta_x$,
then $\Omega(M,x)$ is also flat and $\Omega(M,x)$ is contained in $H_x$ the hyperplane through the origin parallel to $\Delta_x$.
Furthermore,   $\Sigma(M,x)$ and $\Omega(M,x)$ are clearly homothetic, where $\Delta_x\cap  M$ is a symmetric convex body whose boundary is $\Sigma(M,x)$ and its center lies at  $[x,-x]\cap \Delta_x$. 
Moreover, $\Omega(M,x)$  is symmetric with center the origin. 

\medskip

\begin{prop}\label{thmGJMV}
Let $M$ be a symmetric $n$-dimensional convex bodies and let $K$ be a symmetric smooth convex body with $M\subset$ int$K$, $n\geq 3$. Suppose that $M$ looks supersymmetric from $\mbox{bd} K$  and $K$ is almost free with respect to $M$. Then $M$  is an ellipsoid.
\end{prop}

Before giving the proof of Proposition \ref{thmGJMV}, we need several lemmas and establish a notation. Note that by hypothesis, we can conclude that $M$ is strictly convex. So from now, in this section, $M$ is strictly convex body.
We shall first prove the theorem when $n=3$.
\bigskip

 Remember that for every $x\in \mbox{bd}K$, there is a plane $\Delta_x$ such that  $\Delta_x\cap \mbox{bd} M=\Sigma (M,x)$.  
 If $0\in\Delta_x$, that is $\Delta_x=H_x$, then the points of $\Delta_x\cap\mbox{bd} M$ are singular points of the boundary of $M$.
 
 Let 
  $$\eta: \mbox{bd} K\to \mathbb S^2, \quad \mbox{  be the following map: }$$
 $\eta(x)$ is the normal vector to $\Delta_x$ pointing to $x$.

 \begin{lema}\label{homeo}
The map $\eta: \mbox{bd} K\to \mathbb S^2$ is suprayective.
  \end{lema}
 \pf Suppose first $M$ is smooth. Clearly $\eta: \mbox{bd} K\to S^2$ is a continuous map, moreover $\eta$ is injective. This is so because if for $x\not=x' \in \mbox{bd} K$ and $\Delta_x$, $\Delta_{x'}$ are parallel (by smoothness, never the same), then without loss of generality, we may assume that
 $x'$ lies in the interior of $C_xM$,  the cone of M from $x$.  
 Since $M\subset\mbox{int }K$ there is an interior point $q$ of $K$ such that  $x'\in [x,q]$, therefore, $x'$ must be an interior point of $K$, which is a contradiction.  Since $\bd K$ is compact, this implies that,  
  $\eta: \bd K\to \eta(\bd K)$  is a homeomorphism and consequently, $H_{2}(\eta(\bd K), \mathbb Z)=\mathbb Z$. On the other hand, any proper subset of $ \mathbb S^2$ has $2$-dimensional homology equal to zero. Hence  $\eta(\bd K)= \mathbb S^2$, as desired.
  
  In general, when $M$ is not necessarily smooth, the map $\eta: \bd K\to \mathbb S^2$ is not  an homeomorphism, but it is always suprayective. The reason is that 
if $\bd M$ is not smooth, it may happen that $\Delta_x=\Delta_x'$, for $x\not= x'\in \bd K$. If this is the case, 
$\{w \in \bd K\mid \Delta_w=\Delta_x\}$ is contractible, but this implies that the inverse images of the map $\eta$ are contractible. Consequently (see \cite{Lack}), 
$H_{2}(\eta(\bd K), \mathbb Z)=\mathbb Z$ and since any proper subset of $ \mathbb S^2$ has $2$-dimensional homology equal to zero, hence  $\eta(\bd K)= \mathbb S^2$.
\qed

\bigskip

Fix a point $x\in \bd K$  for which $0\notin \Delta_x$ and remember that $\Delta_x \cap \bd M =\Sigma(M,x)$. 
 
 \begin{lema}\label{diana} 
  Suppose $0\notin \Delta_x$.  If $y\in H_x\cap\bd K$,  then
$\Delta_y \mbox{ is parallel to  } [x,-x]$.

\end{lema}

\pf  Let us prove first that  $\Sigma (M,x) \cap \Sigma(M,y)\not=\emptyset.$  If $\Sigma (M,x) \cap \Sigma(M,y)=\emptyset,$  then $\Sigma(M,y)\subset C_xM$ and therefore $y$ lies inside $\Omega(M,x)$, which is a contradiction to the fact that $K$ is almost free with respect to $M$. 
We note that this is the only place in the proof where the almost free hypothesis is used.

Consider now the non-empty chord $I=\Delta_x\cap\Delta_y\cap M$ of $M$. See Figure \ref{para}.
First of all, by symmetry of $M$, and  central symmetry of the sections $\Delta_x\cap K$ and $\Delta_{-x}\cap K$,  there a vector $\tau$ such that 
 $$(\Delta_{-x}\cap K) + \tau  = (\Delta_x\cap K),$$
Furthermore, $\tau$ is parallel to $[x,-x]$ because  in this line segment lie the center of both $\Delta_x\cap K$ and $\Delta_{-x}\cap K$. 
On the other hand, $\Delta_y\cap K$ is also centrally symmetric. Note that the center of $\Delta_y\cap K$ lies in $H_x$ which is at the middle between $\Delta_x$ and $\Delta_{-x}$. Hence, there is a chord $I'$ of $\Delta_y\cap K$ of the same length of $I$. 
Since $0\notin \Delta_x$ and $M$ is strictly convex, $I+\tau= I'$, $\Delta_y$ is parallel to $\tau$ and consequently parallel to $[x,-x]$.
 \qed

\begin{lema}\label{dianalaura}
Suppose $0\notin \Delta_x$ and let $y_0  \in H_x\cap\bd K$.  Then $\Delta_{y_O}\cap H_x$ is an affine line of symmetry of the convex figure $\Delta_{y_O}\cap M$.
\end{lema}

\pf We shall first prove that for every chord $[a,b]$ of $\Delta_{y_O}\cap M$ parallel to $[x,-x]$,  the intersection of the two support lines of the convex figure $\Delta_{y_O}\cap M$ in $\Delta_{y_O}$ at the point $a$ and $b$, lie in the line $\Delta_{y_O}\cap H_x$.  If this is so, 
By Lemma 3 of \cite{JJ} and because $\Delta_{y_O}\cap M$ is centrally symmetric, 
 the reflection of  $\Delta_{y_O}$ in the direction of $[x,-x]$, with mirror line $\Delta_{y_O}\cap H_x$ sends $\Delta_{y_O}\cap M$ to itself as desired. 

\begin{figure}[H]
\centering
\includegraphics [width=1.1\textwidth]{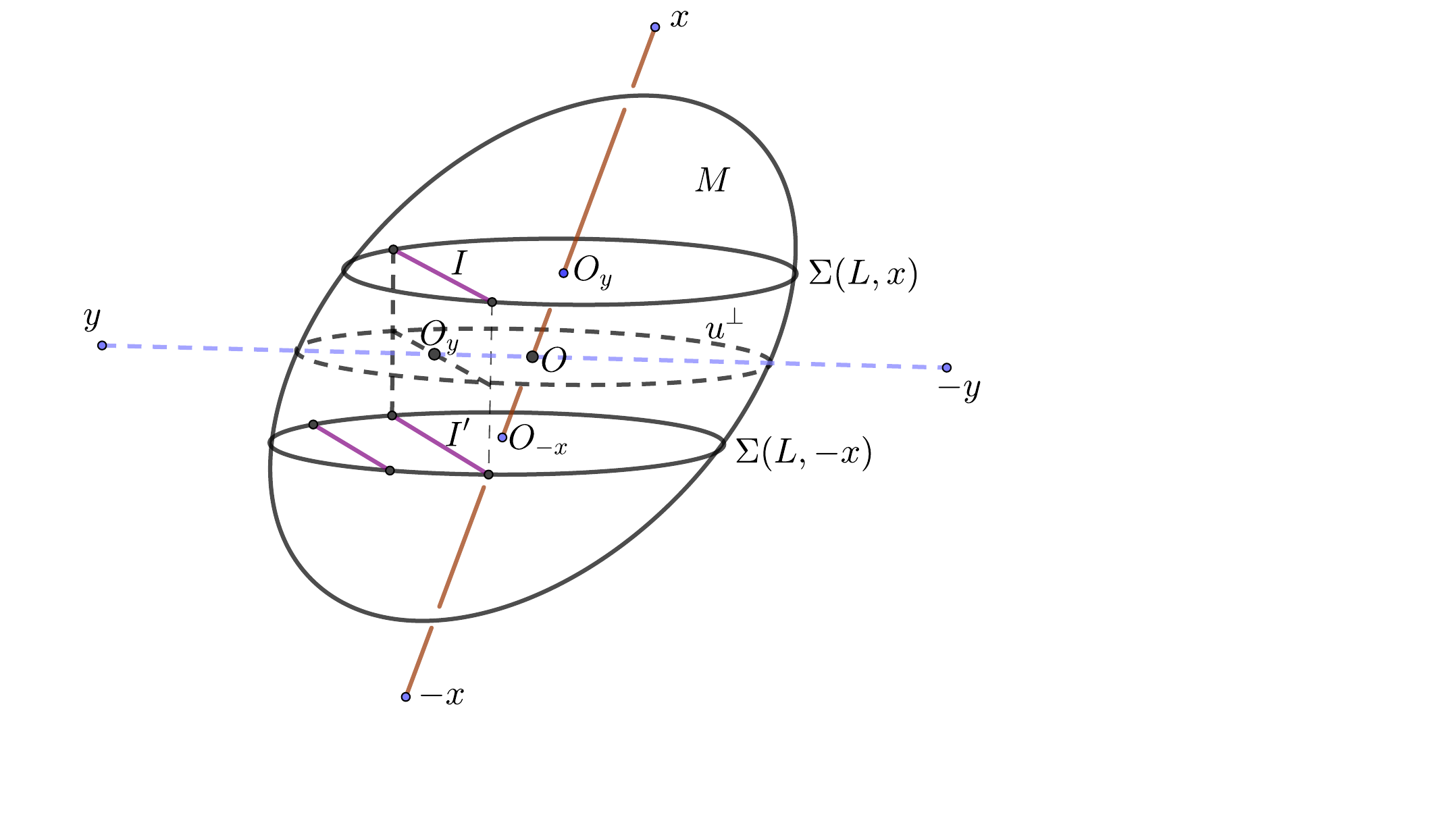}\label{para}
\caption{}\label{para}
\end{figure}

Let $L$ be a line through $y_0$ in $H_x$ such that $L\cap\mbox{int}M=\emptyset$.  If $L$ is not tangent to $K$, then $L\cap\mbox{bd}K=\{y_0,y\}$, where $y_0\not= y\in \mbox{bd}K$.
By Lemma \ref{diana}, $\Delta_y$ is parallel to $[-x,x]$ and therefore $\Sigma(M,L)=\{a,b\}$, where $\Delta_{y_0}\cap\Delta_y\cap M=[a,b]$ and $[ab]$ is parallel to $[-x,x]$.
Hence the plane $\mbox{Aff}\{y_0,y, a\}$ is tangent to $M$ at $a$ and the plane $\mbox{Aff}\{y_0,y, b\}$ is tangent to $M$ at $b$.  Moreover,  $\ell_a:=\mbox{Aff}\{y_0,y, a\}\cap \Delta_{y_0}$ is a line of $\Delta_{y_0}$ tangent to the convex figure
$\Delta_{y_0}\cap M$ at $a$ and $\ell_b:=\mbox{Aff}\{y_0,y, b\}\cap \Delta_{y_0}$ is a line of $\Delta_{y_0}$ tangent to the convex figure
$\Delta_{y_0}\cap M$ at $b$.  Note that $\ell_a\cap\ell_b=L\cap \Delta_{y_0}\in \Delta_{y_0}\cap H_x$, because $L\subset H$.
If $L$ is tangent to $K$ at $y_0$,   $\Sigma(M,L)=\{a,b\}$, where the chord $[a, b]$ of $\Delta{y_0}\cap M$ is parallel to $[-x, x]$. Furthermore, the intersection of the tangent lines of  $\Delta{y_0}\cap M$
at $a$ and $b$ lies in $\Delta_{y_0}\cap H_x$.

Let $[w_0, w_1]=\Delta{y_0}\cap H_x \cap M$ and let $L_0$ be the line through $y_0$ and $w_0$ contained in $H_x$. Similarly, let $L_1$ be the line through $y_0$ and $w_1$ contained in $H_x$. 
If the line $L$ varies from $L_0$ to $L_1$ within the space of lines contained in $H_x$ through $y_0$,  by continuity, the corresponding chord $\{a,b\}=\Sigma(M, L)$ of $\Delta{y_0}\cap M$ parallel to $[-x, x]$, by continuity,  varies between all 
possible chords of $\Delta{y_0}\cap M$ parallel to $[-x, x]$. By Lemma 3 of \cite{JJ} and since $\Delta{y_0}\cap M$ is centrally symmetric, we have that $\Delta_{y_O}\cap H_x$ is an affine line of symmetry of the convex figure $\Delta_{y_O}\cap M$.
\qed

\bigskip
 
\noindent \emph{Proof of Proposition \ref{thmGJMV}}.  Let us first give the proof for $n=3$.  Let $y_0\in \mbox{bd}K$ such that the line through $y_0$ and the origin intersect bd$M$ in a regular point. 
We shall prove that any line $\ell$ through the center of $\Delta_{y_O}\cap M$ is an affine line of symmetry for $\Delta_{y_O}\cap M$. 
Note that the center of $\Delta_{y_O}\cap M$ is 
 $\Delta_{y_0}\cap [0, y_0]$. Let $H$ be the plane through $\ell$ and the origin in such a way that $H$ contains the origin. By Lemma \ref{homeo}, there is $x\in \bd K$ such that $H=H_x$.  
 Furthermore,  $0\notin \Delta_x$, because the line through $y_0$ and the origin intersect bd$M$ in a regular point.
  By Lemma \ref{dianalaura}, $\Delta_{y_O}\cap H_x=\ell$ is an affine line of symmetry of the convex figure $\Delta_{y_O}\cap M$. 
Consequently, for every line $\ell$ through the center of $\Delta_{y_O}\cap M$ in $\Delta_{y_O}$, there a line $\psi(\ell)$  through the center of $\Delta_{y_O}\cap M$ in $\Delta_{y_O}$ such that the midpoint of every chord of 
 $\Delta_{y_O}\cap M$ parallel to $\psi(\ell)$ lies in $\ell$. Note that  $\psi:\mathbb{RP}^1\to\mathbb{RP}^1$ is a well defined continuous map which is also injective. This implies that $\psi(\mathbb{RP}^1)$ is homeomorphic to a circle. 
 Since any proper subset of $\mathbb{RP}^1$ is not a circle, $\psi$ must be a homeomorphism and hence, by Brunn's Theorem 2.12.1 in \cite{MMO},  
$\Delta_{y_O}\cap M$ is an ellipse.  
Finally, the set of points $y_0$ in bd$K$ for which the line through $y_0$ and the origin intersect bd$M$ in a regular point is dense in bd$K$, by continuity, we have that $M$ looks ellipsoidal from  $y_0$,
for every $y_0\in \mbox{bd}K$,   and consequently, by Theorem 2 of \cite{BG},  $M$ is an ellipsoid. 

If $n>3$,  the proof is by induction on $n$. Suppose the theorem holds for $n-1 \geq 3$, we shall prove it for $n$.  
Let $\Gamma$ be a hyperplane through the origin.  We shall prove that for every point  $x\in \bd K \cap \Gamma$,  the graze $\Sigma(\Gamma\cap M, x)$ of $\Gamma\cap M$ in $\Gamma$ is flat. 
For every $x\in \bd K \cap \Gamma$, we know that the grace $\Sigma (M,x)$ is contained in the hyperplane $\Delta_x$.
Indeed, it is easy to see that $x\notin \Delta_x$. Therefore, 
$\Sigma (\Gamma \cap M,x)=\Gamma\cap(M,x)\subset \Gamma\cap \Delta_x$, where of course $\Gamma\cap \Delta_x$ is an 
$(n-2)$-dimensional plane.  Here it is important to note that $\Sigma (\Gamma \cap M,x)$ denotes de graze of $\Gamma \cap M$ in $\Gamma$ from $x$. 
Moreover, since the center of $\Delta_x\cap M$ lies in $[0,x]$, then also the center of $\Gamma\cap\Delta_x\cap M$. 
Similarly, $\Gamma\cap K $ is almost free with respect to $\Gamma \cap M$ because $K$ is almost free with respect to $M$. Consequently, by induction 
$\Gamma\cap M$ is an ellipsoid. 
Since this is true for every hyperplane $\Gamma$ through the origin, we obtain that $M$ is an ellipsoid.
\qed

\medskip

\noindent \emph{Proof of Theorem \ref{super}.}   Let $\phi:\R^n\to \R^n$ be a linear isomorphism that sends $\mathcal E$ to a ball $B_r$.  We may verify that since $M$ looks supersymmetric from every point of the boundary of $\mathcal E$, then 
$\phi(M)$ looks supersymmetric from every point of  $\mathbb S_r$, Furthermore, $\sqrt{2} \phi( M)\subset \mbox{int}B_r$.
Hence, we may assume without loss of generality that $\mathcal E$ is the ball $B_r$.   Consequently, for every point $x\in \mathbb S_r$, the graze $\Sigma(M,x)$ is   contained in the hyperplane 
$\Delta_x$. Let $H_x$ be the plane through the origin parallel to $\Delta_x$.  Remember  that $S(M,x)\cap S(M,-x) =\Omega(M,x) \subset 
H_x$. Let us first prove that $\Omega(M,x) \subset H_x\cap \mathbb B_r$. If this is so,  $\mathbb B_r$ is  almost free with respect to $M$ and hence, by Proposition \ref{thmGJMV}, $M$ is an ellipsoid.  Let $y\in (S(K,x)\cap S(K,-x))$
and consider the parallelepiped with vertices $\{x,y, -x,-y\}$  We first prove that the angle $\measuredangle -yxy$ is less than or equal to
$\pi/2$. This is so because $\sqrt{2}M\subset B_r$. 
On the other hand, $\measuredangle -yxy \leq \pi/2$ implies that $\|y\|\leq \|x\| $. Hence $y\in \mathbb B_r$ and therefore $\Omega(M,x)\subset  B_r$. This proves that $\mathbb B_r$ is  almost free with respect to $M$ as we wished.
\qed

\medskip

The ideas of the proof of Proposition \ref{thmGJMV} come from \cite{GJMV}. In fact, the proof of Lemma 1 in \cite{GJMV} has a gap and therefore the hypothesis was used in the definition of supersymmetry and then Proposition \ref{thmGJMV} was used in the proof of Theorem \ref{super} 

\subsection{ Homothetic ellipsoids}

In \cite{Gru}, P. Gruber showed that among all convex billiard tables in $\mathbb R^n$, $n\geq 3$, only ellipsoids have convex caustics. Caustics of an ellipsoid are precisely the confocal ellipsoids contained in its interior. In his proof, Gruber points out that if $K$ is a billiard table in $\mathbb R^n$ and $M$ is a convex body contained in the interior of $K$, then $M$ is a caustic of $K$ if and only if, for every regular point $p\in \bd K$, the cone $C(M,p)$ is symmetric with respect to the normal line of $K$ at $p$ and therefore the supporting hyperplane of $K$ at $p$ determines the axis of the supporting cone of $M$ at $p$. 
Our next theorem gives a partial answer to Conjecture \ref{conjBG} when we include an additional hypothesis following the spirit of Gruber's Theorem, because  the grazes of $M$ are sections parallel to the support hyperplanes. 

\begin{teo}\label{homo}
Let $M, K$ be convex bodies in $\R^n$, suppose that $M$ is smooth and that $M\subset$ int$K, n\geq3 $. If for every $x\in \mbox{bd} K$, there is a support hyperplane $T_x K$ of $K$ at $x$  and a parallel hyperplane $H$ of $T_x K$ such that
the graze $\Sigma(M,x)= H\cap \mbox{bd M}$, then $M$ and $K$ are homothetic and concentric ellipsoids.
\end{teo}

In the proof of the above result, we will use the following Olevjanischnikoff's Theorem \cite{Ol2} concerning homothetic ellipsoids.

\begin{teo}[Olevjanischnikoff] \label{thmOle}
Let $K$ and $L$ be convex bodies in $\R^n$, with $L\subset$ int$K$, such that every hypersection of $K$ tangent to $L$ is centrally symmetric and its centre belongs to $\mbox{} L$, then $K$ and $L$ are homothetic and concentric ellipsoids.
\end{teo}

\noindent \emph{Proof of Theorem \ref{homo}}.    First of all, note that our hypothesis implies immediately that $M$ is strictly convex. 
Let $x\in \mbox{bd} M$ and let $T_x M$ be the support hyperplane of $M$ at $x$. We shall prove that $T_x M\cap K$ is centrally symmetric with center at $x$. By Theorem \ref{thmOle}, this will conclude our proof.  Let $[p,q]$ be a chord of  
$T_x M\cap K$ through $x$. We first prove that $[p,q]$ is a diametral chord of $T_x M\cap K$.  Since $p,q\in \mbox{bd} K$, by hypothesis, there are hyperplanes $\Gamma_p$ and  $\Gamma_q$, such that $\Gamma_p \cap \mbox{bd} M=\Sigma(M, p)$
and $\Gamma_q \cap \mbox{bd} M=\Sigma(M, q)$. Note first that $\Sigma (M,p)\cap \Sigma (M,q)=\{x\}$. This is so because if $y\not= x $ belongs to  $\Sigma (M,p)\cap \Sigma (M,q)=\{x\}$, then the plane generated be $[p,q]$ and $y$ is tangent to $M$, contradicting that $M$ is strictly convex.
Therefore $\Gamma_p\cap \Sigma(M, q)=\{x\}$. This is so because if $y\in \Gamma_p\cap \Sigma(M, q)$, then $y\in \mbox{bd} M$ and hence by the above $y=x$. Consequently,  $\Gamma_p\cap (\Gamma_q\cap M)=\{x\}$ and therefore the $(n-1)$-plane $\Gamma_p\cap \Gamma_q$ is tangent to 
$M$. On the other hand, the fact that $M$ is smooth implies that $(\Gamma_p\cap \Gamma_q)\subset T_x M$.  Clearly $T_p K\cap T_x M$ is a support $(n-2)$-plane to $T_x M\cap K$ at $p$ and $T_p K\cap T_x M$ is parallel to $\Gamma_p\cap \Gamma_q$, 
because  $T_p K$ is parallel to $\Gamma_p$. Similarly  $T_q K\cap T_x M$ is a support $(n-2)$-plane to $T_x M\cap K$ at $q$ and $T_q K\cap T_x M$ is parallel also to $\Gamma_p\cap \Gamma_q$. Therefore, $T_q K\cap T_x M$ and  $T_p K\cap T_x M$ are parallel support $(n-2)$-planes of 
$T_x M\cap K$ at $q$ and $p$, respectively.  Since this is true for every chord $[p,q]$ of $T_x M\cap K$ through $x$, by Hammer's Theorem cite{H}, $T_x M\cap K$ is centrally symmetric with center at $x$, as we wished. \qed

\bigskip

\noindent {\bf Acknowledgments.} Luis Montejano acknowledge  support  from  PAPIIT-UNAM under project 35-IN112124.

\end{document}